\documentclass[11pt]{amsart}%
\usepackage{amsmath}
\usepackage{amssymb}
\usepackage{amsthm}
\usepackage[mathscr]{eucal}
\usepackage{mathrsfs}
\usepackage{psfrag}
\usepackage{fullpage}
\usepackage{color}
\usepackage{eucal}
\usepackage[dvips]{graphicx}
\usepackage{amsfonts}%
\usepackage{amsaddr}
\providecommand{\U}[1]{\protect\rule{.1in}{.1in}}
\newtheorem{proposition}{Proposition}[section]
\newtheorem{theorem}[proposition]{Theorem}

\newtheorem{definition}[proposition]{Definition}
\newtheorem{remark}[proposition]{Remark}

\newtheorem{condition}[proposition]{Condition}

\numberwithin{equation}{section}
\numberwithin{proposition}{section}

\newcommand{\R}{\mathbb{R}}

\numberwithin{equation}{section}
\numberwithin{proposition}{section}

\begin{document}
\title{Variance reduction for irreversible Langevin samplers and diffusion on graphs}
\author{Luc Rey-Bellet}
\address{Department of  Mathematics and Statistics\\
University of Massachusetts Amherst, Amherst, MA, 01003}
\email{luc@math.umass.edu}

\author{Konstantinos Spiliopoulos}
\address{Department of  Mathematics and Statistics\\
Boston University, Boston, MA, 02215}
\email{kspiliop@math.bu.edu}
\thanks{K.S. was partially supported by the National Science Foundation
(DMS 1312124). LRB was partially supported by the National Science Foundation (DMS 1109316) and the Department of Energy-ASCR
(ER 26161)}

\date{\today}

\maketitle

\begin{abstract}
In recent papers it has been demonstrated that sampling a Gibbs distribution from an appropriate time-irreversible Langevin process is,
from several points of view,  advantageous when compared to sampling from a time-reversible one.
Adding an appropriate irreversible drift to the overdamped Langevin equation results  in a larger large deviations rate function for the empirical measure
of the process, a smaller variance for the long time average of  observables of the process,  as well as a larger spectral gap. In this work, we concentrate on
irreversible Langevin samplers with a drift of increasing intensity. The asymptotic variance  is monotonically decreasing with respect to the growth of the
drift and we characterize its limiting behavior. For a Gibbs measure whose potential has one or
more critical points, adding a large irreversible drift results in a decomposition of the process in a slow and fast component  with fast motion along the level
sets of the potential and slow motion in the orthogonal direction. This result helps understanding the variance reduction, which can be explained  at the
process level by the induced fast motion of the process along the level sets of the potential. The limit of the asymptotic variance as the magnitude of the irreversible perturbation grows
is the  asymptotic variance associated to the limiting slow motion. The latter is a diffusion process on a graph.
\end{abstract}

\section{Introduction}

It is often the case that one is given a high dimensional distribution $\pi(dx)$, which is known only up to normalizing constants, a state space
$E$ and an observable $f$ and the goal is to compute an integral of the form $\bar{f}=\int_{E}f(x)\pi(dx)$. Typically, such integrals cannot be computed in
closed forms, so one has to resort to approximations. What is typically done is to construct a Markov process $X_{t}$, which has $\pi$ as its invariant
distribution. Then under the assumption of positive recurrence, the ergodic theorem guarantees that  for any $f\in L^{1}(\pi)$,
\begin{equation*}
\lim_{t\rightarrow\infty}\frac{1}{t}\int_{0}^{t}f(X_{s})ds=\int_{E}f(x)\pi(dx), \quad\text{a.s.}
\end{equation*}
Then, the estimator $\frac{1}{t}\int_{0}^{t}f(X_{s})ds$ is used to approximate the integral of interest. However, the degree of accuracy of such an approximation depends on both the choice of the Markov process $X_{t}$ and on the criterion used for comparison.

Let us assume that $\pi(dx)$ is of Gibbs type on the state space $E$ and in particular that
\begin{equation}
\pi(dx)=Z^{-1}e^{-\frac{U(x)}{\beta}}dx, \quad{\rm  where }\quad Z=\int_{E}e^{-\frac{U(x)}{\beta}}dx.\label{Eq:GibbsMeasure}
\end{equation}
A Markov process that has $\pi$ as its invariant distribution is the time-reversible Langevin diffusion
\begin{equation}\label{reversibleS}
dX_{t}=-\nabla U(X_{t})dt+\sqrt{2\beta}dW_{t} \,.
\end{equation}
However, as it has been argued in the literature \cite{HwangMaSheu2005,ReyBelletSpiliopoulos2014}, it is advantageous to use appropriate irreversible diffusions of the form
\[
dX_{t}=\left[-\nabla U(X_{t})+ C(X_{t})\right]dt+\sqrt{2\beta}dW_{t}
\]
and if the vector fields $C$ satisfy ${\rm div}( Ce^{-U/\beta})=0$ or equivalently
\[
{\rm div} C = \beta^{-1} C \nabla U \,,
\]
then the measure $\pi$ is still invariant.

The main result of \cite{HwangMaSheu2005} is that under certain conditions the absolute value of the second largest eigenvalue  of the Markov
semigroup in $L^2(\pi)$  decreases when $C\neq0$, which naturally implies faster convergence to equilibrium. In \cite{DupuisDoll1},
the Donsker-Varadhan large deviations rate function  \cite{DonskerVaradhan1975} has been proposed as a natural tool to compare the convergence to
equilibrium for ergodic averages and  it  was used to analyze parallel tempering type algorithms.
In \cite{ReyBelletSpiliopoulos2014}, this criterion is used as a guide to design and analyze
non-reversible Markov processes and compare them with reversible ones. It is proven that the large deviations rate function monotonically increases under
the addition of an irreversible drift. Moreover upon connecting the large deviations rate function with the asymptotic variance of the estimator, it is proven in
\cite{ReyBelletSpiliopoulos2014}, that adding a drift $C$ also decreases the asymptotic variance
of the estimator, in the sense that
\[
\sigma^{2}_{f, C}=\lim_{t\rightarrow\infty} t \text{Var}\left(\frac{1}{t}\int_{0}^{t}f(X_{s})ds\right),
\]
is smaller than the asymptotic variance with $C=0$, i.e, smaller than $\sigma^{2}_{f, 0}$.


The goal of the current work is to study the situation when the additional drift has the form  $\frac{1}{\epsilon}C(x)$ and to consider what happens when $\epsilon \to 0$.
A  similar question but from a  different perspective has been studied in \cite{ConstantinKiselevRyshikZlatos2008,FrankeHwangPaiSheu2010}. There, the
authors  studied the behavior of the spectral gap for diffusions on compact manifolds  with $U=0$ and a one-parameter families of perturbations
$\frac{1}{\epsilon} C$ for some divergence free vector field $C$.  In those papers the behavior of the spectral gap is related to the
ergodic properties of the flow generated by $C$ (for example if the flow is weak-mixing then the second largest  eigenvalue tends to $0$ as $\epsilon \to 0$).
In the present paper we want to understand the effect that increasing $1/\epsilon$ has on the asymptotic variance and  on  the paths of  $X_{t}$.

We find that the asymptotic variance of the estimator is monotonically decreasing in $\delta=1/\epsilon$. Using an averaging argument,
we characterize the limiting asymptotic variance as $\epsilon\downarrow 0$. Focusing on the case where the potential $U(x)$ has one 
or more critical points, the irreversible perturbations with a small $\epsilon$ induce a fast motion on the constant potential surface and 
slow motion in the orthogonal direction. Using the theory of diffusions on  graphs and the related averaging principle as developed in 
a series of works \cite{BrinFreidlin2000, FreidlinWeber2004,FreidlinWentzell1993, FW2}, we identify the limiting motion of the slow 
component. The fast motion on constant potential surfaces decreases the variance as the phase space is explored faster. The limit of 
the asymptotic variance as $\epsilon\to 0$ is the asymptotic variance of a one-dimensional estimation problem on a graph, which is where 
the slow component of the process lives in the limit as $\epsilon\downarrow 0$.
%
%

Upon completion of this work, we became aware of the recent paper \cite{HwangNormandWu2014} where an alternative expression of the limiting asymptotic
variance is provided, see also Remark \ref{R:resultOfHwang}. The methods of \cite{HwangNormandWu2014} are analytical and the characterization of the limiting variance is expressed as the
projection to a kernel of a certain operator. Our approach provides complementary information and  is dynamical since we are using an averaging principle which allows us to make direct
connections with the limiting behavior of the underlying process itself, relating the limiting variance with an estimation problem on an one-dimensional graph.

The rest of the paper is organized as follows. In Section \ref{S:MainResults} we formulate the problem precisely  and present our main results.
The averaging problem treating the limit of the slow component of the process is discussed in  Section \ref{S:AveragingProblem}.  Due to the special structure
of  the model, one can perform explicit computations and thus derive precise results. The formula for the limit of the asymptotic variance as the perturbation
grows is in Section \ref{S:VarianceBehavior}.
Numerical simulations illustrating the theoretical findings are in Section \ref{S:Simulations}.

\section{Statement of the problem and main results}\label{S:MainResults}
The papers \cite{ConstantinKiselevRyshikZlatos2008,FrankeHwangPaiSheu2010,HwangMaSheu2005,ReyBelletSpiliopoulos2014} motivate to look at  a one parameter family of irreversible drifts $\frac{1}{\epsilon}C(x)$. For this purpose, we consider the model
\begin{equation}
dX^{\epsilon}_{t}=\left[-\nabla U(X^{\epsilon}_{t})+ \frac{1}{\epsilon} C(X^{\epsilon}_{t})\right]dt+\sqrt{2\beta}dW_{t}\label{Eq:IrreversibleDiffusion}
\end{equation}
where $\frac{1}{\epsilon} C$ is a one-parameter family of vector fields,  $\epsilon\in\mathbb{R}$ and the vector field $C$ satisfies
${\rm div}( Ce^{-U/\beta})=0$.   As  mentioned in the introduction, the invariant measure is maintained if the vector fields $C$  satisfies ${\rm div}( Ce^{-U/\beta})=0$, or equivalently
\[
{\rm div} C = \beta^{-1} C \nabla U \,.
\]

A convenient choice, which we assume henceforth, is to pick $C$ such that
\[
{\rm div} C =0 \,, \quad {\rm and} \quad  C \nabla U =0 \,.
\]

This is not the most general choice for $C$, but it has the advantage that allows to choose $C$ independently of $\beta$. A standard choice of $C(x)$ is $C(x)=S\nabla U(x)$, where $S$ is any antisymmetric matrix. A  more
elaborate discussion on other possible choices of $C(x)$ can be found in \cite{ReyBelletSpiliopoulos2014}. The meaning of these conditions
is straightforward: the flow generated by $C$ must preserve Lebesgue measure since it it is divergence-fee but since $U$ is a constant of the motion,  the
micro-canonical measure  on the surfaces $\{U = z\}$ are preserved as well. 

We assume here that the diffusion process  $X^{\epsilon}_{t}$ is on a $d$-dimensional compact smooth manifold without boundary and 
that $U$ and $C$ are sufficiently smooth.  The ergodicity of $X^\epsilon_t$ implies that  the empirical measure
\begin{equation*}
\pi^{\epsilon}_{t} \equiv \frac{1}{t}\int_{0}^{t} \delta _{X^{\epsilon}_s} \, ds,
\end{equation*}
converges to $\pi$ almost surely as $t\rightarrow\infty$.  Under our assumptions we have a large deviation principle (uniformly in the initial condition)
for the family of measures $\pi^{\epsilon}_t$, which we write, symbolically, as
\begin{equation*}
\mathbb{P} \left\{ \pi^{\epsilon}_{t}  \approx \mu \right\} \asymp \exp \left(- t I_{1/\epsilon }(\mu)\right)
\end{equation*}
where $\asymp$ denotes logarithmic equivalence. (Since $C$ is fixed we have suppressed the dependence of $I_{1/{\epsilon}}$ on $C$.)
The rate function $I_{1/{\epsilon}}(\mu)$
quantifies the exponential rate at which the (random) measure $\pi^{\epsilon}_t$ converges to $\pi$ as $t\rightarrow\infty$.
It is proven in  \cite{ReyBelletSpiliopoulos2014}  (see Theorem \ref{T:measure2}  below) that the rate function $I_{ 1/\epsilon }(\mu)$ is quadratic in
$\epsilon$.

The information in $I_{1/\epsilon}(\mu)$ can be used to study the rate of convergence of observables. If $f \in \mathcal{C}(E; \mathbb{R})$
then we have the large deviation principle
\begin{equation*}
\mathbb{P} \left\{ \frac{1}{t}\int_0^t f(X^{\epsilon}_s) \, ds  \approx \ell \right\} \asymp \exp \left( - t \tilde{I}_{f,1/\epsilon}(\ell) \right) \,,
\end{equation*}
where, by the contraction principle, we have
\begin{equation*}
\tilde{I}_{f, 1/\epsilon}(\ell)=\inf_{\mu\in\mathcal{P}(E)}\left\{I_{1/\epsilon }(\mu) \, ; \int f d\mu=\ell\right\} \,.
\end{equation*}

An alternative formula for $\tilde{I}_{f, 1/\epsilon}(\ell)$ is in terms of the Legendre transform of the maximal eigenvalue $\lambda(\gamma)$ of the Feynmann-Kac
semigroup $T^t_\gamma h(x) ={\bf E}_x \left[ e^ {\gamma \int_0^t f(X^\epsilon_s) \, ds} h( X^\epsilon_t)\right]$ acting on ${\mathcal C}(E;\mathbb{R})$.
For  $\ell$ in the range of $f$ we have
\begin{equation}
\tilde{I}_{f, 1/\epsilon}(\ell) = \sup_{\gamma} \left[\gamma \ell - \lambda(\gamma) \right] \,=\,  \hat{\gamma}(\ell) \ell - \lambda( \hat{\gamma}(\ell))  \quad {\rm where}  \quad
\lambda'(\hat{\gamma})= \ell  \,.
\end{equation}

For $f \in L^2(\pi)$ with $\bar{f} = \int f d\pi$ the asymptotic variance is given by, see e.g. Proposition I.V.1.3 in \cite{AsmussenGlynn2007},
\[
\sigma^{2}_{f,
1/\epsilon}=\lim_{t\rightarrow\infty} t \text{Var}\left(\frac{1}{t}\int_{0}^{t}f(X^\epsilon_{s})ds\right) \,=\, 2\int_{0}^{\infty} \mathbb{E}_{\pi}\left[(f(X^{\epsilon}_0)- \bar{f})(f(X^{\epsilon}_{t}) - \bar{f})\right]dt
\]
and if  $ f \in  {\mathcal C}(E;\mathbb{R})$ it is related to the rate function $\tilde{I}_{f, 1/\epsilon}(\ell)$  by
\[
\tilde{I}''_{f, 1/\epsilon}( \bar{f}) \,=\, \frac{1}{2 \sigma^{2}_{f, 1/\epsilon}} \,.
\]

We have the following theorem from \cite{ReyBelletSpiliopoulos2014}.
\begin{theorem}[Theorem 2.3, 2.4, and 2.6 of \cite{ReyBelletSpiliopoulos2014}]\label{T:measure2} Assume that $E$ is a smooth connected manifold without
boundary and assume that  for some $\alpha >0$ , $U \in {\mathcal C}^{(2+\alpha)}(E)$ and $C \in {\mathcal C}^{(1+\alpha)}(E)$ are such that
${\rm div}(C e^{-U/\beta}) =0$. Then $\pi$ is the invariant measure of the process $X^{\epsilon}_{t}$ and the following hold.
\begin{enumerate}
\item Let $\mu(dx) =p(x)dx$ be a measure with positive density $p \in \mathcal{C}^{(2 + \alpha)}(E)$. Then we have
\begin{equation*}
I_{1/\epsilon }(\mu) = I_0(\mu) + \frac{1}{\epsilon^{2}} K(\mu)  \,,
\end{equation*}
where the functional $K(\mu)$ is positive and strictly positive if and only if $\text{div}\left(p(x)C(x)\right)\not=0$.  It takes the explicit form
\begin{equation*}
K(\mu)=\frac{1}{2}\int_{E}\left|\nabla \xi(x)\right|^{2}d\mu(x) \,,
\end{equation*}
where $\xi$ is the unique solution (up to a constant) of the equation
${\textrm div}\left[p\left(C+\nabla \xi\right)\right]=0$.
\item  We have $\tilde{I}_{f,\frac{1}{\epsilon}C }(\ell) \ge \tilde{I}_{f,0}(\ell)$ and generically the inequality is strict:
for  $f \in {\mathcal C}^{(\alpha)}(E)$ we have
$\tilde{I}_{f,\frac{1}{\epsilon}C }(\ell) = \tilde{I}_{f,0}(\ell)$ if and only if there exists $G$ invariant under the vector field $C$ such that
\begin{equation*}
\widehat{\gamma}(\ell) f \,=\, \mathcal{H}(G+U)\,,\label{Eq:NonlinearPoissonEquation}
\end{equation*}
where $\mathcal{H}(G+U)=e^{-(G+U)}\mathcal{L}_{0}e^{G+U}$ with $\mathcal{L}_{0}=\beta \Delta -\nabla U\nabla$.
\item For the asymptotic variance we have
\[
\sigma^{2}_{f, 1/\epsilon} \leq \sigma^{2}_{f, 0}
\]
with strict inequality if  $\tilde{I}_{f,\frac{1}{\epsilon}C }(\ell)>\tilde{I}_{f,0}(\ell)$ in a neighborhood of $\bar{f}=\int_{E}f(x)\pi(dx)$, but excluding $\bar{f}$.
More generally the map  $|\epsilon|\mapsto \sigma^{2}_{f, \frac{1}{\epsilon}}$ is a monotone increasing function and thus
its limit as $\epsilon\to 0$ exists.
\end{enumerate}
\end{theorem}

For notational convenience we shall write,  from now on, 
\[
 \sigma^{2}_{f}(\epsilon)=\sigma^{2}_{f,\frac{1}{\epsilon}},
\]
and, without loss of generality, we may and will assume that $\bar{f}=0$.
We want to characterize the behavior of the asymptotic variance as $\epsilon\downarrow 0$, namely to find
$\lim_{\epsilon\downarrow 0}\sigma^{2}_{f}(\epsilon)$, and to connect it with the limiting behavior of the trajectories
$X^{\epsilon}_{t}$ as $\epsilon\downarrow 0$.

Notice that we can write
\begin{align*}
\sigma^{2}_{f}(\epsilon)&=2\int_{0}^{\infty}\mathbb{E}_{\pi}\left[f(X^{\epsilon}_{0})f(X^{\epsilon}_{t})\right]dt=
2\int_{E}\Phi^{\epsilon}(x)f(x) \pi(dx)
\end{align*}
where $\Phi^{\epsilon}(x)$ is the unique solution (up to constants) of the Poisson equation
\begin{equation}
-\mathcal{L}_{\epsilon}\Phi^{\epsilon}(x)=f(x)\label{Eq:Poisson}
\end{equation}
with $\int_{E}\Phi^{\epsilon}(x)\pi(dx)=0$. Here $\mathcal{L}_{\epsilon}$ is the infinitesimal generator of the process $X^{\epsilon}_{t}$ given by (\ref{Eq:IrreversibleDiffusion}). Dissipativity of the operator $\mathcal{L}_{\epsilon}$ implies that
\begin{align*}
\sigma^{2}_{f}(\epsilon)&=2\int_{E}\Phi^{\epsilon}(x)f(x) \pi(dx)=2\int_{E}\Phi^{\epsilon}(x)(-\mathcal{L}_{\epsilon}\Phi^{\epsilon}(x)) \pi(dx)\geq 0.
\end{align*}

It is easy to see that as $\epsilon\rightarrow 0$, the solution  to (\ref{Eq:Poisson}) becomes constant on the stream lines of $C(x)$. In particular, if we multiply  (\ref{Eq:Poisson})
by $\epsilon$ and then take $\epsilon\rightarrow 0$, we formally obtain that $C(x)\cdot\nabla \Phi(x)=0$. Making this rigorous is the purpose of Sections \ref{S:AveragingProblem} and \ref{S:VarianceBehavior}.
Note that this heuristic argument makes it immediately clear that for $f\in\mathcal{L}^{2}(\pi)$
\begin{equation*}
 \text{if } \left[\text{Ker}(C\cdot\nabla)=\{0\}\right] \Longrightarrow [\sigma^{2}_{f}(0)=0].
\end{equation*}

Hence, we next investigate what happens in the non-trivial case, i.e., when $\text{Ker}(C\cdot\nabla)\neq\{0\}$. Our goal is to relate the behavior of the variance as $\epsilon\rightarrow 0$
with that of the process $X^{\epsilon}_{t}$. It turns out that as $\epsilon\rightarrow 0$, the behavior of the process $X^{\epsilon}_{t}$ can be decomposed into fast motion on
 the constant
potential surface and slow motion in the orthogonal direction. Using averaging principle, see for example   \cite{BrinFreidlin2000,FreidlinWeber2004,FreidlinWentzell1993,  FW2}, we can identify the limiting behavior of the slow component  and then compute
$\sigma^{2}_{f}(0)=\lim_{\epsilon\rightarrow0}\sigma^{2}_{f}(\epsilon)$. To be more precise, the slow component of the motion can be characterized in the limit as $\epsilon\rightarrow 0$ as an
one-dimensional
Markov process on a graph and $\sigma^{2}_{f}(0)$ turns out to be associated with the asymptotic variance of an estimation problem on the graph itself.

The precise behavior of the process $X^{\epsilon}_{t}$ as $\epsilon\rightarrow 0$ has been investigated in \cite{FreidlinWeber2004, FW2}.
We will state the precise assumptions and results
in Section \ref{S:AveragingProblem}. However, let us briefly review the result here.
Following \cite{FW2}, we consider a finite graph $\Gamma$, which represents the structure of the level sets of
the potential function $U$ on $E$. To construct the graph we identify the points that belong to the connected components of each of the level sets of $U$. We assume that $U$ has finitely many
non-degenerate critical points and that each connected level set component of $U$ contains at most one critical point. In that way, each of the domains that is bounded by the separatrices
gets mapped into an edge of the graph. At the same time the separatrices gets mapped to the vertexes of $\Gamma$. In particular, exterior vertexes correspond to minima of $U$, whereas interior
vertexes correspond to saddle points of $U$. Edges of $\Gamma$ are indexed by $I_{1},\cdots, I_{m}$. Each point on $\Gamma$ is indexed by a pair
 $y=(z,i)$ where $z$ is the value of $U$ on the level set corresponding to $y$ and $i$ is the edge number containing $y$. Clearly the pair $y=(z,i)$  forms a global coordinate on $\Gamma$. Let $Q:E\mapsto\Gamma$ with $Q(x)=(U(x),i(x))$ be the corresponding projection on the graph.

Consider the process $X^{\epsilon}_{t}$ on $E$ to be the solution to (\ref{Eq:RegularizedIrreversibleLangevin}). As we describe in Section \ref{S:AveragingProblem}, the additional
noise that appears in (\ref{Eq:RegularizedIrreversibleLangevin}) is added for regularization reasons and the limit does not depend on it. Moreover, the additional regularizing
noise is not needed if $C(x)$ generates an ergodic dynamical system with a unique invariant measure within each connected component 
of the level sets of $U$.  It is proven in Theorem 2.1 of \cite{FreidlinWeber2004} and in Theorem 2.2 of \cite{FW2} (if $C$ is an Hamiltonian 
vector field) that the process $Q(X^{\epsilon}_{t})$ converges to a certain Markov process on $\Gamma$ with continuous trajectories, which is exponentially mixing. The precise result will be stated
in Section \ref{S:AveragingProblem}, but roughly speaking it goes as follows.

\begin{theorem}[\cite{FreidlinWeber2004, FW2}] \label{T:AveragingResult1}
 Let us assume that Conditions \ref{A:Assumption2}, \ref{A:Assumption3} and \ref{A:Assumption4} hold. Then, for any $0<T<\infty$, the process $Y^{\epsilon}_{t}=Q(X^{\epsilon}_{t})$, where $X^{\epsilon}_{t}$ satisfies (\ref{Eq:RegularizedIrreversibleLangevin}), converges weakly in $\mathcal{C}([0,T],\Gamma)$ to a Markov process,
denoted by $Y_{t}$, on $\Gamma$ with continuous trajectories, which is exponentially mixing.
\end{theorem}

We remark here that the classical Freidlin-Wentzell theory, see \cite{FW2}, assumes that the diffusion on $\mathbb{R}^d$ with 
$\lim_{|x|\to\infty} U(x)=\infty$. But the results apply in the compact case as well.

Based on Theorem \ref{T:AveragingResult1}, we can then establish the limiting behavior of the asymptotic variance as $\epsilon\downarrow 0$ and then make connections to estimation problems on the graph.
In particular we have the following result which is discussed and proven in Section \ref{S:VarianceBehavior}.
\begin{theorem}\label{T:LimitingVariance}
Let us assume that Conditions \ref{A:Assumption2}, \ref{A:Assumption3} and \ref{A:Assumption4} hold and let $Y_{t}$ be the continuous Markov process on the graph $\Gamma$ indicated in Theorem  \ref{T:LimitOnGraphb}.
Let
$f\in \mathcal{C}^{2+\alpha}(E)
$ such that $\bar{f}=0$. For $(z,i)\in\Gamma$, define $\widehat{f}(z,i)$ to be the average of $f$
on the graph $\Gamma$ over the corresponding
connected component of the level set $U$ (see equation (\ref{Eq:AverageOnGraph}) for precise definition). Then, we have that
$\sigma^{2}_{f}(0)=\lim_{\epsilon\rightarrow0}\sigma^{2}_{f}(\epsilon)$, where
\begin{equation}
 \sigma^{2}_{f}(0)=2\int_{0}^{\infty}\mathbb{E}_{\mu}\left[\widehat{f}(Y_{0})\widehat{f}(Y_{t})\right]dt\label{Eq:LimitingVariance}
\end{equation}
and $\mu=\pi\circ\Gamma^{-1}$ is the invariant measure of the process $Y$ on $\Gamma$.
\end{theorem}

Theorem \ref{T:LimitingVariance} is proven in Section \ref{S:VarianceBehavior}. It is straightforward to see that this is the asymptotic variance of an ergodic average on the graph. In particular, we have
\begin{equation}
 \sigma^{2}_{f}(0)=\lim_{t\rightarrow\infty} t \text{Var}\left(\frac{1}{t}\int_{0}^{t}\widehat{f}(Y_{s})ds\right).\label{Eq:LimitingVarianceOnTheGraph}
\end{equation}

\begin{remark}\label{R:resultOfHwang}
For completeness purposes, we briefly recall here the result of \cite{HwangNormandWu2014} that is related to the present situation. When the irreversible perturbation is of the form $\frac{1}{\epsilon}C$ with $C$ chosen such that ${\rm div}( Ce^{-U/\beta})=0$, Theorem 4.3 of \cite{HwangNormandWu2014} states that $\sigma^{2}_{f}(0)=2 P\mathcal{L}_{0}^{-1/2}f$. Here $\mathcal{L}_{0}$ is the infinitesimal generator of the process $X_{t}$ in (\ref{reversibleS}) and $P$ is the projection on $\text{Ker}\left( i\mathcal{L}^{-1/2}_{0}\left(C\cdot\nabla\right)\mathcal{L}^{-1/2}_{0}\right)$.

The methods of \cite{HwangNormandWu2014} are analytical based on an analysis of the spectrum of the operator and of the related spectral measure.  The methods of our paper are dynamic, formula (\ref{Eq:LimitingVariance}) is derived using an averaging principle and it is valid under the constraint ${\rm div} C =C \nabla U =0$. Since the methodologies for deriving the two results are very different, one naturally obtains different equivalent expressions for the limiting asymptotic variance.  It is of great interest to understand how one can go from one formulation to the other, at least in the case ${\rm div} C =C \nabla U =0$. Doing so, would then also allow one to connect the projection operators that appear in \cite{HwangNormandWu2014} with objects such as diffusions processes on graphs.
\end{remark}

\section{The averaging problem}\label{S:AveragingProblem}
In this section we discuss the behavior of the process $X^{\epsilon}_{t}$ as $\epsilon\downarrow 0$. Such problems have been studied in \cite{BrinFreidlin2000,FreidlinWentzell1993, FreidlinWeber2004, FW2}
and we recall here the results which are relevant  to us. It turns out that as $\epsilon\rightarrow 0$, the behavior of the process $X^{\epsilon}_{t}$ can be decomposed into fast motion on
the constant potential surface and slow motion in the orthogonal direction. 

Let us consider the level set 
\begin{equation}
d(z)=\left\{x\in E: U(x)=z\right\}
\end{equation}
and denote by $d_{i}(z)$  the connected components of $d(z)$, i.e.,
\begin{equation}
d(z)=\bigcup_{i}d_{i}(z)
\end{equation}

Then, we let  $\Gamma$ be the graph which is homeomorphic to the set of connected components $d_{i}(z)$ of the level sets $d(z)$. 
Exterior vertexes correspond to minima of $U$, whereas interior vertexes correspond to saddle points of $U$. The edges of $\Gamma$ 
are indexed by $I_{1},\cdots, I_{m}$. Each point on $\Gamma$ is indexed by a pair $y=(z,i)$ where $z$ is the value of $U$ on the level 
set corresponding to $y$ and $i$ is the edge number containing $y$. Clearly the pair $y=(z,i)$
forms a global coordinate on $\Gamma$. Let $Q: E \mapsto\Gamma$ with $Q(x)=(U(x),i(x))$
be the corresponding projection on the graph. For an edge $I_{k}$ and a vertex $O_{j}$ we write $I_{k}\sim O_{j}$ if $O_{j}$
lies at the boundary of the edge $I_{k}$. We endow the tree $\Gamma$ with the natural topology. It is known that $\Gamma$ forms a 
graph  with interior vertexes of order  two or three, see for example \cite{FW2}.

Let us next consider $X^{\epsilon}_{t}$ with an additional artificial noise component in the fast dynamics, i.e.,
\begin{equation}\label{Eq:RegularizedIrreversibleLangevin}
dX^{\epsilon}_{t}=\left[-\nabla U(X^{\epsilon}_{t})dt+\sqrt{2\beta}dW_{t}\right]+ \left[ \frac{1}{\epsilon}  \tilde{C}(X^{\epsilon}_{t})dt+\sqrt{\frac{\kappa}{\epsilon}} \sigma(X^{\epsilon}_{t})dW^{o}_{t}\right]
\end{equation}
where $W$ and $W^{o}$ are independent standard Wiener processes, and we have defined
\[
\tilde{C}(x)=C(x)+\frac{\kappa}{2}\sum_{j=1}^{d}\frac{\partial \left[\sigma\sigma^{T}(x)\right]_{j,i}}{\partial x_{j}}.
\]
If $\kappa=0$ then we get the process $X^{\epsilon}_{t}$ that we have been considering until now.

We make several technical assumptions on $C(x)$, $U(x)$ and $\sigma(x)$ in order to guarantee that the averaging principle 
applies. We make these assumptions in order to guarantee that the fast process has a unique invariant measure and will have $U$ 
as a smooth first integral. If $\kappa=0$ then the fast motion is the deterministic dynamical system $\dot{x}_{t}=C(x_{t})$ and 
$X^{\epsilon}_{t}$ is  a random perturbation of this dynamical system. For example, if $d$ is even we can take $C$ to be the
Hamiltonian vector field $C(x)=\bar{\nabla}U(x)$. If $\kappa>0$ we have random perturbations of diffusion processes with a 
conservation law.  
In order to guarantee the existence of a unique invariant measure for the fast dynamics we assume:
\begin{condition}\label{A:Assumption2}
In dimension $d=2$, we take $\kappa\geq 0$. In dimension $d>2$, we either assume that the dynamical system $\dot{x}_{t}=C(x_{t})$ 
has a unique invariant measure on each connected component $d_{i}(z)$, in which case $\kappa\geq 0$, or otherwise we assume 
that $\kappa>0$. 
\end{condition}

As far as the potential function $U(x)$ and the perturbation $C(x)$ are concerned, we shall assume:
\begin{condition}\label{A:Assumption3}
\begin{enumerate}
\item{There exists $a>0$ such that $U\in\mathcal{C}^{(2+a)}(E)$ and $C\in \mathcal{C}^{(1+a)}(E)$.}
\item{$div C(x)=0$ and $C(x)\cdot\nabla U(x)=0$.}
\item{$U$ has a finite number of critical points $x_{1},\cdots, x_{m}$ and at these points the Hessian matrix is non-degenerate.}
\item{There is at most one critical point for each connected level set component of $U$.}
\item{If $x_{k}$ is a critical point of $U$, then there exists $d_{k}>0$ such that $C(x)\leq d_{k}|x-x_{k}|$.}
\item{If $d=2$ and $\kappa=0$, then $C(x)=0$ implies $\nabla U(x)=0$ and for any saddle point $x_{k}$ of $U(x)$, there exists a constant $c_{k}>0$ such that $|C(x)|\geq c_{k}|x-x_{k}|$.}
\end{enumerate}
\end{condition}

In regards to the additional artificial perturbation by the noise $W^{o}_{t}$, i.e., when $\kappa>0$, we assume:
\begin{condition}\label{A:Assumption4}
\begin{enumerate}
\item{The matrix $\sigma(x)\sigma^{T}(x)$ is non-negative definite, symmetric with smooth entries.}
\item{$\sigma(x)\sigma^{T}(x)\nabla U(x)=0$ for all $x\in E$.}
\item{For any $x\in E$ such that $\xi\cdot\nabla U(x)=0$ we have that $\lambda_{1}(x)|\xi|^{2}\leq \xi^{T}\sigma(x)\sigma^{T}(x)\xi\leq \lambda_{2}(x)|\xi|^{2}$ where $\lambda_{1}(x)>0$ if $\nabla U(x)\neq 0$ and there exists a constant $K$ such that $\lambda_{2}(x)<K$ for all $x\in E$.
Moreover if $x_{k}$ is a critical point for $U$, then there are positive constants $k_{1},k_{2}$ such that for all $x$ in a neighborhood of $x_{k}$
    \begin{equation*}
    \lambda_{1}(x)\geq k_{1} |x-x_{k}|^{2}, \text{ and }     \lambda_{2}(x)\leq k_{2} |x-x_{k}|^{2}.
    \end{equation*}}
\item{Let $\lambda_{i,k}$ be the eigenvalues of the Hessian of $U(x)$ at the critical points $x_{k}$ where $k=1,\cdots,m$ and $i=1,\cdots, d$. Then we assume that
$\kappa<\left(K\max_{i,k}\lambda_{i,k}\right)^{-1}$.}
\end{enumerate}
\end{condition}

Obviously Condition \ref{A:Assumption4} is relevant in the case $d>2$ and if $C(x)$ does not generate an ergodic dynamical system (since otherwise we can just take $\kappa=0$).
In the  case $\kappa>0$, the procedure of  incorporating an appropriate artificial noise in the system, allows to single out
the correct averaging that should be done in the system. We remark here that the end result does not depend on the additional
regularizing noise, since $\sigma(x)$ does not appear in the limiting dynamics.

It is clear that the dynamics can be decomposed in a fast component and a slow component. The fast motion corresponds to the infinitesimal generator
\begin{equation*}
\hat{\mathcal{L}}g(x)=\tilde{C}(x)\nabla g(x)+\frac{\kappa}{2}\text{tr}\left[\sigma\sigma^{T}(x)\nabla^{2}g(x)\right] \,.
\end{equation*}

Let us write $\hat{X}_{t}$ for the diffusion process that has infinitesimal generator $\hat{\mathcal{L}}$. Conditions \ref{A:Assumption3} and \ref{A:Assumption4} guarantee that with probability one, if the initial point of $\hat{X}$ is in a connected component $d_{i}(z)$, then $\hat{X}_{t}\in d_{i}(z)$ for all $t\geq 0$. Indeed, by It\^{o} formula we have
\begin{equation*}
U(\hat{X}_{t})=U(\hat{X}_{0})+\int_{0}^{t}\hat{\mathcal{L}}U(\hat{X}_{s})ds+\int_{0}^{t}\nabla U(\hat{X}_{s})\sigma(\hat{X}_{s})dW_{s}.
\end{equation*}

Since $C(x)\nabla U(x)=0$ and $\sigma(x)\sigma^{T}(x)\nabla U(x)=0$ we obtain with probability one $\int_{0}^{t}\hat{\mathcal{L}}(\hat{X}_{s})ds=0$. The quadratic variation of the stochastic integral is also zero, due to $\sigma(x)\sigma^{T}(x)\nabla U(x)=0$, which implies that with probability one  $\int_{0}^{t}\nabla U(\hat{X}_{s})\sigma(\hat{X}_{s})dW_{s}=0$. Thus, we indeed get that for all $t\geq 0$ $\hat{X}_{t}\in d_{i}(z)$ given that the initial point belongs to the particular connected component $d_{i}(z)$.

Let us turn now our attention to invariant measures.  Let $m(x)$ be a smooth invariant density with respect to Lebesgue measure for the 
process $\hat{X}_{t}$. Then, the proof of Lemma 2.3 of \cite{FreidlinWeber2004} and the fact that $t\geq 0$ $\hat{X}_{t}\in d_{i}(z)$ if $
\hat{X}_{0}\in d_{i}(z)$ imply that if $(z,i)\in\Gamma$ is not a vertex, there exists a unique invariant measure $\mu_{z,i}$ concentrated on 
the connected component $d_{i}(z)$ of $d(z)$ which takes the form
\begin{equation}
\mu_{z,i}(A)=\frac{1}{T_{i}(z)}\oint_{A}\frac{m(x)}{\left|\nabla U(x) \right|}\ell(dx)
\end{equation}
where $T_{i}(z)=\oint_{d_{i}(z)}\frac{m(x)}{\left|\nabla U(x) \right|}\ell(dx)$. 
Notice that if $(z,i)\in\Gamma$ is not a vertex, then the invariant density on $d_{i}(z)$ is
\begin{equation}
m_{z,i}(x)
=\frac{m(x)}{T_{i}(z)\left|\nabla U(x) \right|},\quad x\in d_{i}(z).
\end{equation}

We remark here that in the case $\kappa>0$, it is relatively easy to see that independently of the form of the matrix $\sigma(x)\sigma^{T}(x)$, the fact that $\text{div}(C)=0$ implies that the Lebesgue measure is invariant for the diffusion process corresponding to the operator  
$\hat{\mathcal{L}}$. Hence, in that case any constant function is an invariant density. Also, in the case $d=2$ and $\kappa=0$, one immediately obtains from Condition \ref{A:Assumption3} that $m(x)=\frac{|\nabla U(x)|}{|C(x)|}$, see Proposition 2.1 in \cite{FreidlinWeber2004}.

Given a sufficiently smooth  function $f(x)$, define its average over the related connected component of the level set of $U(x)$ by
\begin{equation}
\widehat{f}(z,i)=\oint_{d_{i}(z)}f(x)m_{z,i}(x)\ell(dx)=\frac{1}{T_{i}(z)}\oint_{d_{i}(z)}\frac{f(x)}{\left|\nabla U(x) \right|}m(x)\ell(dx)\label{Eq:AverageOnGraph}
\end{equation}

Let us consider the process $Q(X^{\epsilon}_{t})=\left(U(X^{\epsilon}_{t}), i(X^{\epsilon}_{t})\right)$ and consider its limiting behavior.  We write $\mathcal{L}_{0}$ for the infinitesimal generator of the process $X_{t}$ given by (\ref{reversibleS}).
Let us set
\begin{align*}
\widehat{\mathcal{L}_{0}U}(z,i)&=\oint_{d_{i}(z)}\mathcal{L}_{0}U(x)m_{z,i}(x)\ell(dx)=\frac{1}{T_{i}(z)}\oint_{d_{i}(z)}\frac{\mathcal{L}_{0}U(x)}{\left|\nabla U(x)\right|}m(x)\ell(dx),\nonumber\\
\widehat{A}(z,i)&=\oint_{d_{i}(z)}2\beta |\nabla U(x)|^{2}T_{i}(z)m_{z,i}(x)\ell(dx)=\oint_{d_{i}(z)}\frac{2\beta \nabla U(x)\cdot\nabla U(x)}{\left|\nabla U(x)\right|}m(x)\ell(dx)
\end{align*}
and then consider the one-dimensional process $Y_{t}$ which within the branch $I_{i}$ is governed by the infinitesimal generator
\begin{equation*}
\mathcal{L}^{Y}_{i}g(z)= \widehat{\mathcal{L}_{0}U}(z,i) g'(z)+\frac{1}{2}\frac{\widehat{A}(z,i)}{T_{i}(z)}g''(z)
\end{equation*}

Within each edge $I_{i}$ of $\Gamma$, $Q(X^{\epsilon}_{t})$ converges as $\epsilon\downarrow 0$ to a process
with infinitesimal generator $\mathcal{L}^{Y}_{i}$. In order to uniquely define the limiting process, we need to specify the behavior at the vertexes of the tree, which amounts to imposing restrictions on the domain of definition of the generator, say $\mathcal{L}^{Y}$, of the Markov process.

\begin{definition}\label{Def:ProcessOnTree}
We say that $g$ belongs in the domain of definition of $\mathcal{L}^{Y}$, denoted by $\mathcal{D}$, of the diffusion $Y_{\cdot}$, if
\begin{enumerate}
\item{The function $g(z)$ is twice continuously differentiable in the interior of an edge $I_{i}$.}
\item{The function $z\mapsto\mathcal{L}^{Y}_{i}g(z)$ is continuous on $\Gamma$.}
\item{At each interior vertex $O_{j}$ with edges $I_{k}$ that meet at $O_{j}$, the following gluing condition holds
    \begin{equation*}
    \sum_{k: I_{k}\sim O_{j}}\pm b_{jk}D_{k}g(O_{j})=0
    \end{equation*}
    where, if $\gamma_{jk}$ is the separatrices  curves that meet at $O_{j}$,  we have set
    \begin{equation*}
    b_{jk}=\oint_{\gamma_{jk}}\frac{2\beta \left|\nabla U(x)\right|^{2}}{\left|\nabla U(x)\right|}m(x)\ell(dx)=2\beta\oint_{\gamma_{jk}}\left|\nabla U(x)\right|m(x)\ell(dx).
    \end{equation*}
    Here one chooses $+$ or $-$ depending on whether the value of $U$ increases or decreases respectively along
 the edge $I_{k}$ as we approach $O_{j}$. Moreover  $D_{k}$ represents the derivative in the direction of the edge $I_{k}$.
    }
\end{enumerate}
Moreover, within each edge $I_{i}$ the  process $Y_{t}$ is a diffusion process with infinitesimal generator $\mathcal{L}^{Y}_{i}$.
\end{definition}

Consider now the process $Y_{t}$ that has the aforementioned $\mathcal{L}^{Y}$ as its infinitesimal generator with domain of definition $\mathcal{D}$, as defined in Definition \ref{Def:ProcessOnTree}.
 Such a process is a continuous strong Markov process, e.g., Chapter 8 of \cite{FW2}.

Then, for any $T>0$, $Q(X^{\epsilon}_{t})$ converges weakly in $\mathcal{C}([0,T];\Gamma)$ to the process $Y_{t}$ as $\epsilon\downarrow 0$.
In particular, we have the following theorem.
\begin{theorem}[Theorem 2.1 of \cite{FreidlinWeber2004}]\label{T:LimitOnGraphb}
Let $X^{\epsilon}_{t}$ be the process that satisfies (\ref{Eq:RegularizedIrreversibleLangevin}). Assume Conditions \ref{A:Assumption2}, \ref{A:Assumption3} and \ref{A:Assumption4}. Let $T>0$ and consider
 the Markov process $\left\{Y_{t},t\in[0,T]\right\}$ as defined in Definition \ref{Def:ProcessOnTree}. We have
\begin{equation}
 Q(X^{\epsilon}_{\cdot})\rightarrow Y_{\cdot}, \text{ weakly in } \mathcal{C}([0,T];\Gamma), \text{ as }\epsilon\downarrow 0.
\end{equation}
\end{theorem}

It is important to note that the limiting process $Y_{t}$ does not depend on $\sigma(x)$.  Next, we show that in our case of interest and for every $i=1,\cdots,m$, the operator $\mathcal{L}^{Y}_{i}$ which governs the motion of the limiting process within the $I_{i}$ branch of the graph takes a more explicit form.
Let us denote $G_{i}(z)=\text{int}(d_{i}(z))$. Notice that by Gauss theorem we have

\begin{align*}
\widehat{\mathcal{L}_{0}U}(z,i)&=\frac{1}{T_{i}(z)}\oint_{d_{i}(z)}\frac{\mathcal{L}_{0}U(x)}{\left|\bar{\nabla} U(x)\right|}m(x)\ell(dx)=\frac{1}{T_{i}(z)}\oint_{d_{i}(z)}\frac{-\left|\nabla U(x)\right|^{2}+\beta\Delta U(x)}{\left|\nabla U(x)\right|}m(x)\ell(dx)\nonumber\\
&=\frac{1}{T_{i}(z)}\left[-\int_{G_{i}(z)}\text{div}\left(m(x)\nabla U(x) \right)dx+ \beta\oint_{d_{i}(z)}\frac{\Delta U(x)}{\left|\nabla U(x)\right|}m(x)\ell(dx)\right]
\end{align*}
and similarly
\begin{align*}
\widehat{A}_{i}(z)&=\oint_{d_{i}(z)}\frac{2\beta \nabla U(x)\cdot\nabla U(x)}{\left|\bar{\nabla} U(x)\right|}m(x)\ell(dx)=2\beta\int_{G_{i}(z)}\text{div}\left(m(x)\nabla U(x) \right) dx
\end{align*}

Next, we notice that
\begin{align*}
\frac{d}{dz}\int_{G_{i}(z)}\Delta U(x) m(x)dx&=\oint_{d_{i}(z)}\frac{\Delta U(x)}{\left|\nabla U(x)\right|} m(x)\ell(dx).
\end{align*}

Hence, we can write
\begin{align*}
\widehat{\mathcal{L}_{0}U}(z)
&=\frac{1}{T_{i}(z)}\left[-\int_{G_{i}(z)}\mathcal{L}^{*}_{0}m(x)dx+ \frac{1}{2} \widehat{A}_{i}^{\prime}(z)\right]
\end{align*}
where $\mathcal{L}^{*}_{0}$ is the formal adjoint operator to $\mathcal{L}_{0}$. Thus, the infinitesimal generator  $\mathcal{L}^{Y}_{i}$, can be written equivalently  as
\begin{align*}
\mathcal{L}^{Y}_{i}g(z)&=\frac{1}{T_{i}(z)}\left[-\int_{G_{i}(z)}\mathcal{L}^{*}_{0}m(x)dx\right] g'(z)+\frac{1}{2T_{i}(z)}\frac{d}{dz}\left\{\widehat{A}_{i}(z)g'(z)\right\}.
\end{align*}

In particular, if the Lebesgue measure is the invariant measure (e.g., in the case $\kappa>0$), then $m(x)$ is a constant. By denoting $M_{i}(z)=\int_{G_{i}(z)}\Delta U(x)dx>0$,  we can then rewrite $\mathcal{L}^{Y}_{i}$ in the more explicit form
\begin{align*}
\mathcal{L}^{Y}_{i}g(z)&=-\frac{1}{T_{i}(z)}M_{i}(z) g'(z)+\frac{\beta}{T_{i}(z)}\frac{d}{dz}\left\{M_{i}(z) g'(z)\right\}\nonumber\\
&=\frac{-M_{i}(z)+\beta M_{i}^{\prime}(z)}{T_{i}(z)} g'(z)+\beta \frac{M(z)}{T_{i}(z)} g''(z).
\end{align*}


\section{Limiting behavior of the asymptotic variance}\label{S:VarianceBehavior}
In this section we identify $\sigma^{2}_{f}(0)=\lim_{\epsilon\downarrow 0}\sigma^{2}_{f}(\epsilon)$ using the averaging results of Section \ref{S:AveragingProblem}. In particular,
we prove the representation of Theorem \ref{T:LimitingVariance}.

 Recall that without loss of generality we assume $\bar{f}=0$. Our starting point is  the formula
\begin{equation*}
\sigma^{2}_{f}(\epsilon)=2\int_{0}^{\infty}\mathbb{E}_{\pi}\left[f(X^{\epsilon}_{0})f(X^{\epsilon}_{s})\right]ds,
\end{equation*}
where the process $X^{\epsilon}_{t}$ is the unique strong solution of (\ref{Eq:RegularizedIrreversibleLangevin}). Standard PDE arguments, e.g., Section 3.2 of \cite{FreidlinWeber2004}, show that for any point $(z,i)\in I_{i}$ that is not a vertex, the PDE
\begin{equation}
-\hat{\mathcal{L}}u(x)=f(x)-\widehat{f}(U(x),i(x)), \quad\text{for }x\in d_{i}(z)\label{Eq:AuxilliaryPDE}
\end{equation}
has a unique solution (up to constants) $\mathcal{C}^{2+\alpha'}$ solution with $\alpha'\in(0,\alpha)$. We fix the free constant by setting $\widehat{u}(z,i)=0$. Then, the solution $u(x)$ can be written as
\[
u(x)=\int_{0}^{\infty}\mathbb{E}_{x}\left[f(\hat{X}_{s})-\widehat{f}(U(\hat{X}_{s}),i(\hat{X}_{s}))\right]ds.
\]

Moreover, there exist constants $\lambda=\lambda(z,i)>0$ such that for $x\in d_{i}(z)$,
\begin{equation}
|u(x)|\leq \frac{2}{\lambda}\sup_{x\in d_{i}(z)}\left|f(x)-\widehat{f}(U(x),i(x))\right|.\label{Eq:BoundAuxilliaryPDE}
\end{equation}

Let us consider $\theta>0$ small and for an edge $I_{i}$ of the graph set
\[
 I^{\theta}_{i}=\{(z,i)\in I_{i}: \text{dist}((z,i),\partial I_{i})>\theta\}
\]
and define
\[
 \tau_{i}=\min\{t>0: Q(X^{\epsilon}_{t})\notin I^{\theta}_{i})\}
\]

Then, by applying It\^{o} formula to the solution of (\ref{Eq:AuxilliaryPDE}) with stochastic process $X^{\epsilon}_{t}$ one immediately gets that for any $T<\infty$,
for any initial point $x$ that does not belong to any of the separatrices of $U$ and for every $I_{i}$
\begin{equation}\label{Eq:AVeragingLimitWithinBranch}
\lim_{\epsilon\downarrow 0} \sup_{t\in[0,T]}\mathbb{E}_{x}\left[\int_{0}^{t\wedge\tau_{i}}\left(f(X^{\epsilon}_{s})-\widehat{f}(U(X^{\epsilon}_{s}),i(X^{\epsilon}_{s}))\right)ds\right]=0
\end{equation}
uniformly in $x\in D_{i}(z)=\{x\in\R^{d}: Q(x)\subset\text{int} (I_{i})\}$. Indeed, applying It\^{o} formula to $u(x)$ with $x=X^{\epsilon}_{s}$, we have
\begin{align*}
u(X^{\epsilon}_{t})&=u(x_{0})+\int_{0}^{t}\left[-\nabla U(X^{\epsilon}_{s})\nabla u(X^{\epsilon}_{s})+\frac{1}{\epsilon}\hat{\mathcal{L}}u(X^{\epsilon}_{s})
+\beta\Delta
u(X^{\epsilon}_{s})\right]ds\nonumber\\
&+\int_{0}^{t}\sqrt{2\beta}\nabla u(X^{\epsilon}_{s})dW_{s}+\sqrt{\frac{\kappa}{\epsilon}}\int_{0}^{t}\nabla u(X^{\epsilon}_{s})\sigma(X^{\epsilon}_{s})dW^{o}_{s}
\end{align*}
 and recalling that $u(x)$ is solution to (\ref{Eq:AuxilliaryPDE}) we get
\begin{align*}
\int_{0}^{t\wedge\tau_{i}}\left(f(X^{\epsilon}_{s})-\widehat{f}(U(X^{\epsilon}_{s}),i(X^{\epsilon}_{s}))\right)ds &=
\epsilon\left[u(x_{0})-u(X^{\epsilon}_{t\wedge\tau_{i}})\right]+\epsilon\int_{0}^{t\wedge\tau_{i}}\left[-\nabla U\nabla u+\beta\Delta
u\right](X^{\epsilon}_{s})ds\nonumber\\
&\quad+\epsilon\int_{0}^{t\wedge\tau_{i}}\sqrt{2\beta}\nabla u(X^{\epsilon}_{s})dW_{s}+\sqrt{\kappa\epsilon}\int_{0}^{t\wedge\tau_{i}}\nabla u(X^{\epsilon}_{s})\sigma(X^{\epsilon}_{s})dW^{o}_{s}
\end{align*}

Taking expected value, the right hand side of this inequality goes to zero as $\epsilon\downarrow 0$, by (\ref{Eq:BoundAuxilliaryPDE}) and because continuity of the integrands implies that
 Riemann integrals are bounded. Hence,  the averaging result (\ref{Eq:AVeragingLimitWithinBranch}) follows immediately.

 At the same time, by the results of \cite{FreidlinWeber2004, FW2}
the limiting process $Y_{t}$ spends time of Lebesgue measure zero at the interior and exterior vertexes.
For $\zeta>0$ and for a vertex of the graph $O_{j}$, let us define
\[
 D_{j}(\pm\zeta)=\left\{x\in\mathbb{R}^{d}: U(O_{j})-\zeta<U(x)<U(O_{j})+\zeta\right\}
\]

If $O_{j}$ is an exterior vertex
of $\Gamma$, then for every $\eta>0$, there exists $\zeta>0$ such that for sufficiently small $\epsilon$ and for all $x\in D_{j}(\pm\zeta)$ we have that (Lemma 3.6 in \cite{FreidlinWeber2004})
\[
 \mathbb{E}_{x}\tau^{\epsilon}_{j}(\pm\zeta)<\eta
\]
where $\tau^{\epsilon}_{j}(\pm\zeta)$ is the first exit time of $X^{\epsilon}_{t}$ from $D_{j}(\pm\zeta)$. The behavior for an interior vertex is similar. Lemma 3.7 of \cite{FreidlinWeber2004} implies that
if $O_{j}$ is an interior vertex, then for every $\eta>0$, and for all sufficiently small $\zeta>0$
\[
 \mathbb{E}_{x}\tau_{j}^{\epsilon}(\pm\zeta)<\eta\zeta
\]
for sufficiently small $\epsilon$ and for all $x\in D_{j}(\pm\zeta)$.
Recall now that  $Q(x)=(U(x),i(x))$. Then, (\ref{Eq:AVeragingLimitWithinBranch}) and the fact that the process $Y_{t}$ spends time of Lebesgue measure zero at all vertexes,
 imply that for any $t<\infty$
\begin{align}
 \lim_{\epsilon\downarrow 0}\int_{0}^{t}\mathbb{E}_{\pi}\left[f(X^{\epsilon}_{0})f(X^{\epsilon}_{s})\right]ds&=\lim_{\epsilon\downarrow 0}\int_{0}^{t}\mathbb{E}_{\pi}\left[f(X^{\epsilon}_{0})\widehat{f}(U(X^{\epsilon}_{s}),i(X^{\epsilon}_{s})\right]ds\nonumber\\
&=\lim_{\epsilon\downarrow 0}\int_{0}^{t}\mathbb{E}_{\pi}\left[\widehat{f}(Q(X^{\epsilon}_{0}))\widehat{f}(Q(X^{\epsilon}_{s}))\right]ds
\end{align}

Since, the invariant measure of $X^{\epsilon}_{t}$ is the Gibbs measure $\pi$, we have that
the invariant measure of $Y^{\epsilon}_{t}=Q(X^{\epsilon}_{t})$ and of $Y_{t}$ on the graph is the projection of the Gibbs measure $\pi$ on $\Gamma$.
Denoting this invariant measure by $\mu$, we have that for any Borel set $\gamma\subset \Gamma$,  $\mu(\gamma)=\pi(\Gamma^{-1}(\gamma))$. Thus,
by the weak convergence of Theorem \ref{T:LimitOnGraphb} we have that for any $t<\infty$
\begin{align}
 \lim_{\epsilon\downarrow 0}\int_{0}^{t}\mathbb{E}_{\pi}\left[f(X^{\epsilon}_{0})f(X^{\epsilon}_{s})\right]ds&=
\int_{0}^{t}\mathbb{E}_{\mu}\left[\widehat{f}(Y_{0})\widehat{f}(Y_{s})\right]ds.
\end{align}

The strong Markov processes $X^{\epsilon}_{t}$ and $Y_{t}$,  on $E$ and $\Gamma$ respectively,   are
uniform mixing.
This implies that, by selecting $t>0$ to be large enough, we can make the integrals
\[
\int_{t}^{\infty}\mathbb{E}_{\pi}\left[f(X^{\epsilon}_{0})f(X^{\epsilon}_{s})\right]ds\quad \text{and}\quad \int_{t}^{\infty}\mathbb{E}_{\mu}\left[\widehat{f}(Y_{0})\widehat{f}(Y_{s})\right]ds
\]
arbitrarily small. Therefore, we indeed obtain that
\[
 \lim_{\epsilon\rightarrow 0}\sigma^{2}_{f}(\epsilon)=\sigma^{2}_{f}(0)=2\int_{0}^{\infty}\mathbb{E}_{\mu}\left[\widehat{f}(Y_{0})\widehat{f}(Y_{s})\right]ds,
\]
concluding the proof of Theorem \ref{T:LimitingVariance}.

\section{Numerical Simulations}\label{S:Simulations}

In this section, we explore numerically the behavior of the process under  growing perturbations of the drift. In the examples below we fix $\beta=0.1$.

\begin{figure}[t!]
\begin{center}
\includegraphics[scale=0.4]{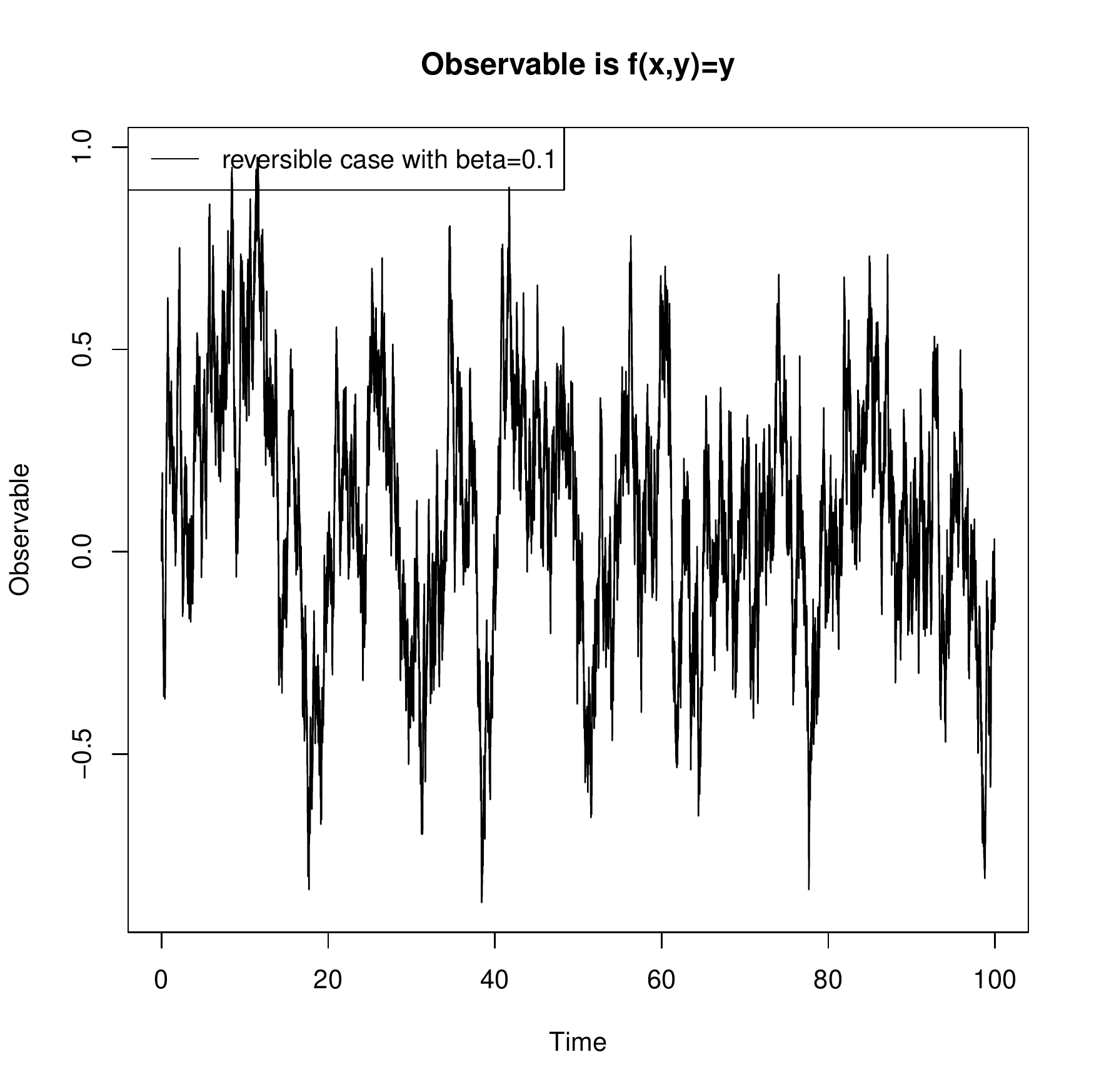}
\includegraphics[scale=0.4]{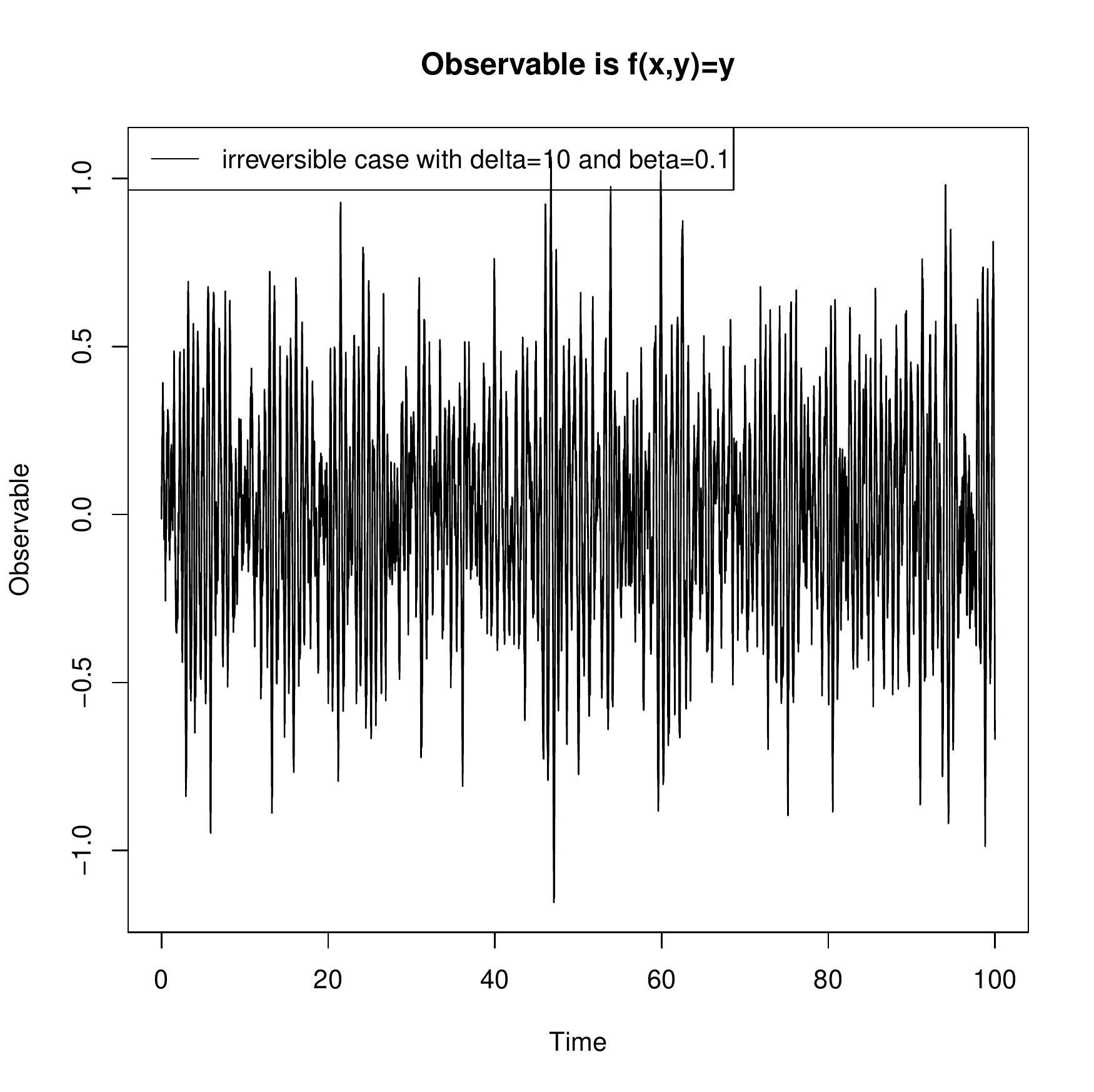}
\end{center}
\caption{\label{Fig1} On the left: reversible case, i.e. $\delta=0$. On the right: irreversible case with $\delta=10$. In both cases $\beta=0.1$.
 }
\end{figure}

The first example that we study is a simple 2-dimensional example where the potential $U$ has a single critical point. In particular,
we define $U(x,y)=\frac{1}{2}x^{2}+\frac{1}{2}y^{2}$.
Let $C(x,y)=S \nabla U(x,y)$, where $S$ is the standard $2\times2$ antisymmetric matrix, i.e., $S_{12}=1$ and $S_{21}=-1$.
In Figure \ref{Fig1}, we see that the more irreversibility one adds (in the sense of increasing the $\delta=1/\epsilon$
parameter in the perturbation $\delta C(x,y)$), the faster the process explores the phase space. Since, it is perhaps convenient to think in terms of $1/\epsilon$ and not $\epsilon$,
we set $\delta=1/\epsilon$.


Furthermore, notice that what the theory predicts is also shown in the numerical simulations. Namely, we observe fast motion long the level sets of the potential and slow motion in the orthogonal direction.

Next, we consider a  two-dimensional SDE with two asymptotically stable points. In particular, we define the potential
$U(x,y)=\frac{1}{4}(x^{2}-1)^{2}+\frac{1}{2}y^{2}$ and we consider growing perturbation $\delta C(x,y)$ with $C(x,y)=S \nabla U(x,y)$.
 Here, $\delta\in \mathbb{R}$ and $S$ is the standard $2\times2$ antisymmetric matrix.

\begin{figure}[t!]
\begin{center}
\includegraphics[scale=0.4]{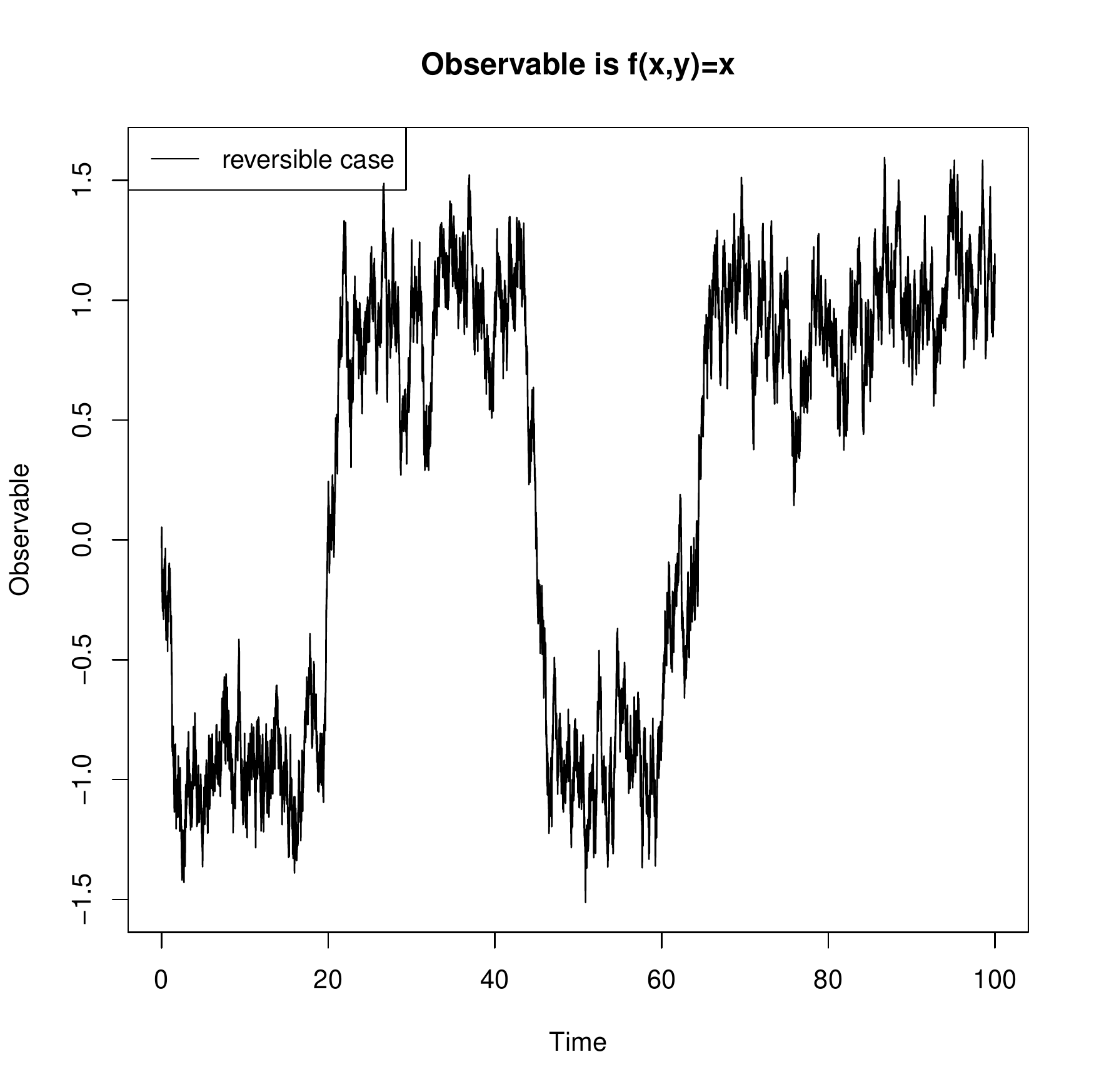}
\includegraphics[scale=0.4]{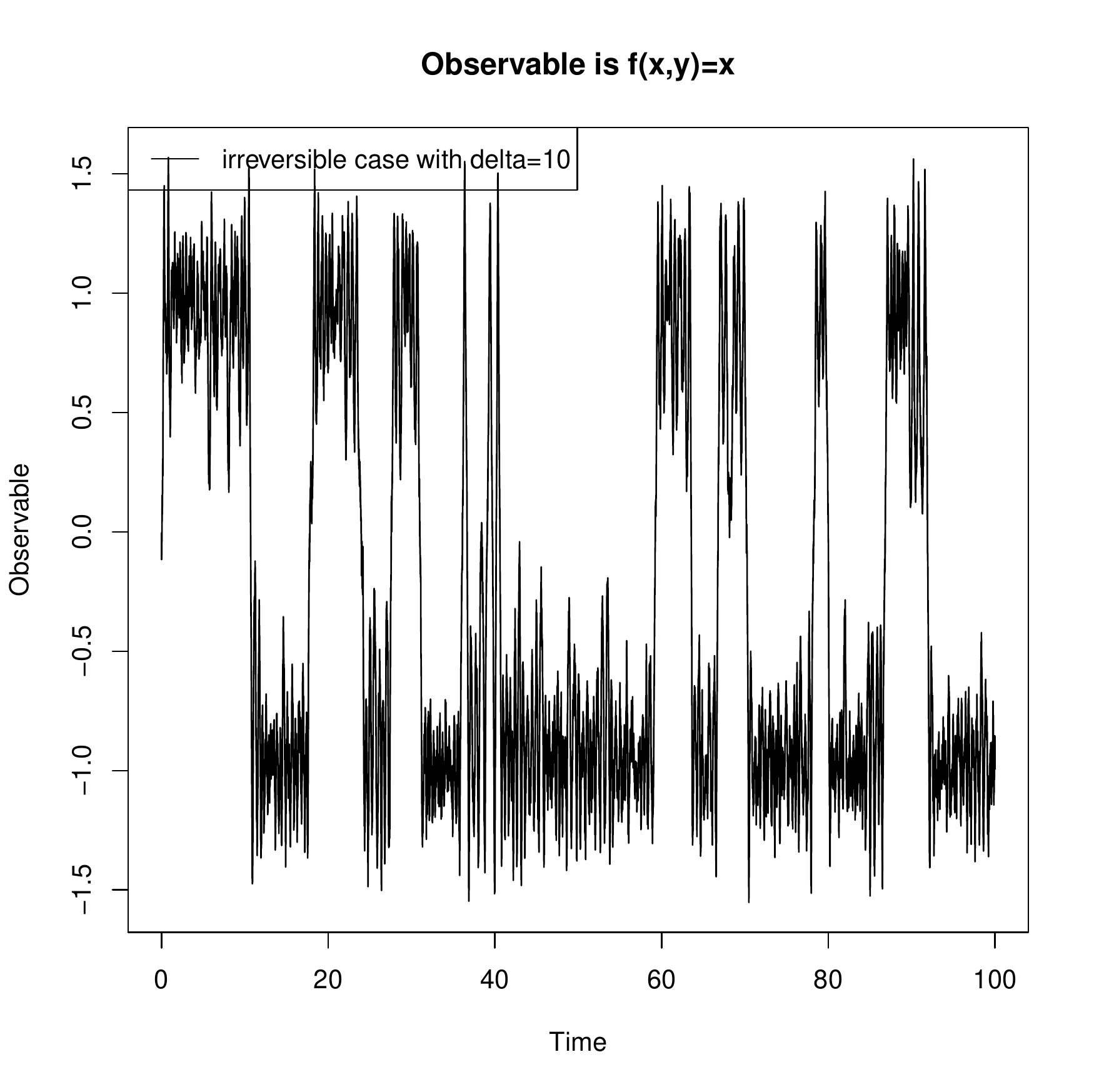}
\end{center}
\caption{\label{Fig7} On the left: reversible case, i.e. $\delta=0$. On the right: irreversible case with $\delta=10$.
 }
\end{figure}

\begin{figure}[t!]
\begin{center}
\includegraphics[scale=0.4]{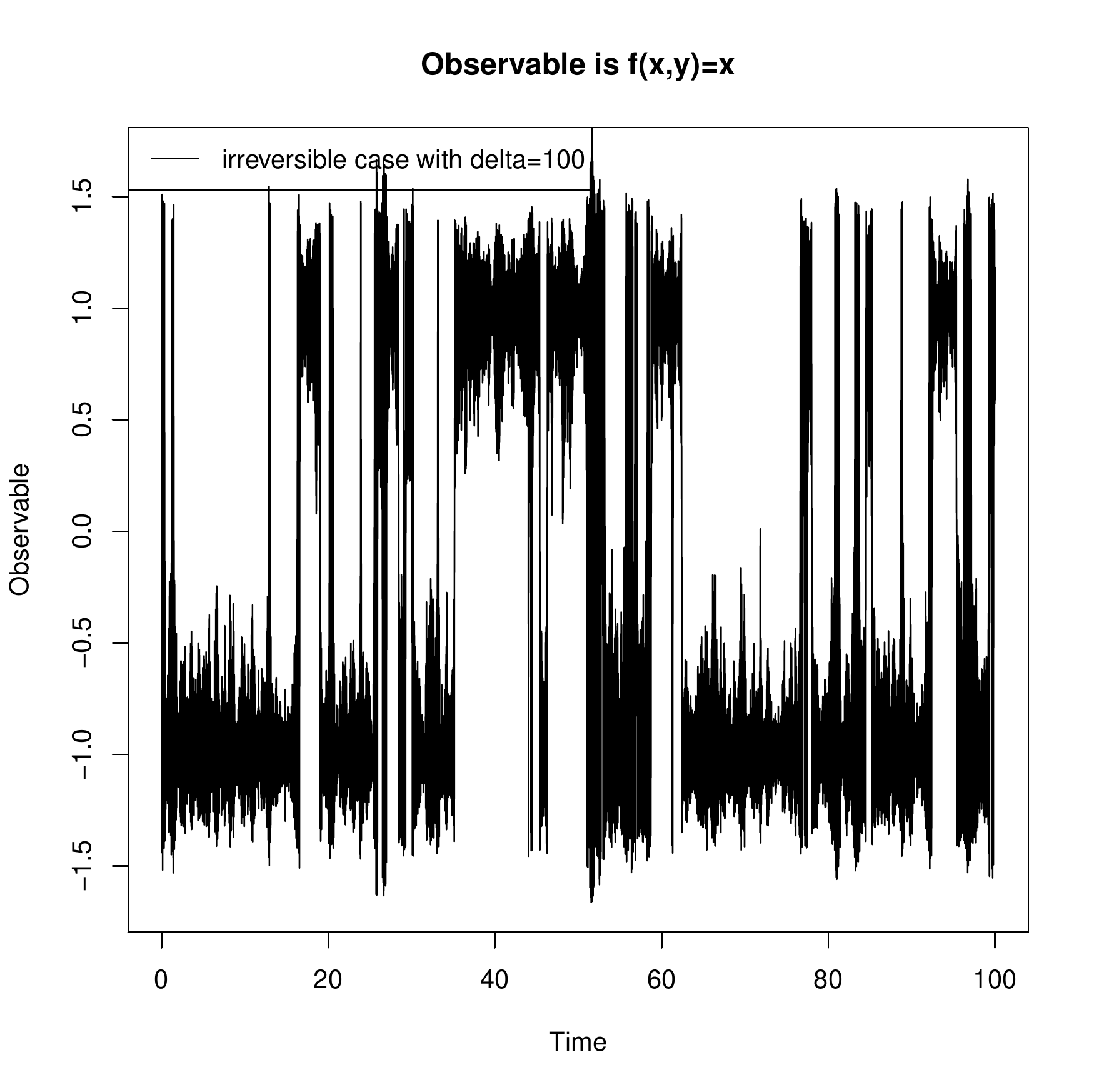}
\includegraphics[scale=0.4]{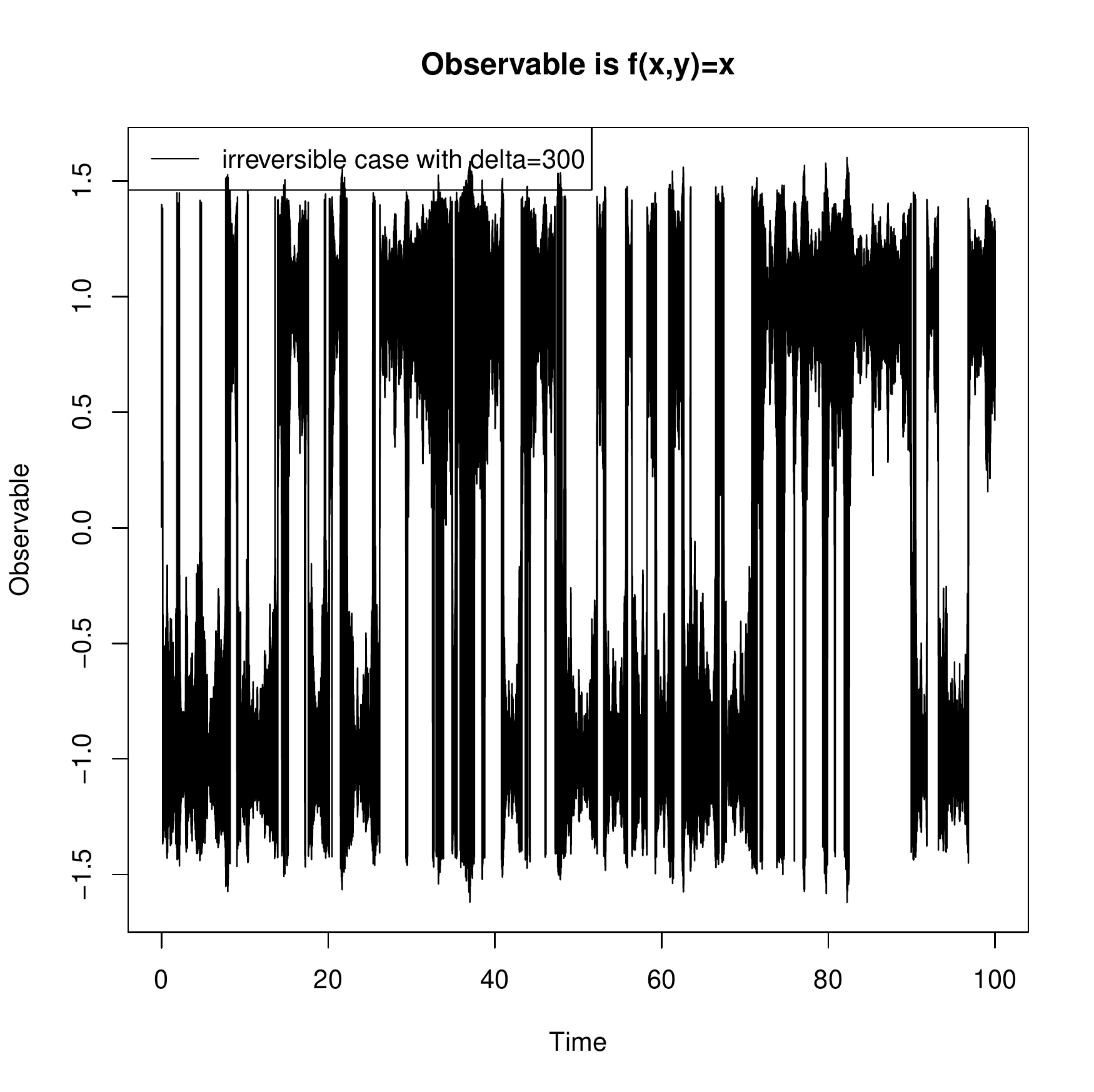}
\end{center}
\caption{\label{Fig8} On the left: irreversible case with $\delta=100$. On the right: irreversible case with $\delta=300$.
 }
\end{figure}

Notice that the potential function $U(x,y)$ has two local minima in $(-1,0)$ and $(1,0)$ and a local maximum at $(0,0)$. In Figures \ref{Fig7} and \ref{Fig8}, we plot the $x-$component of the $(x,y)$ trajectory versus time. We see that the
more irreversibility one adds (in the sense of increasing the $\delta$ parameter in the perturbation $\delta C(x,y)$), the faster the process moves  along the level sets of the
 potential. As Figures \ref{Fig7} and \ref{Fig8} show, in the present case, this means  faster
 switches between the two metastable states.

\begin{figure}[ptb]
\includegraphics[height=4cm,width=9.25cm]{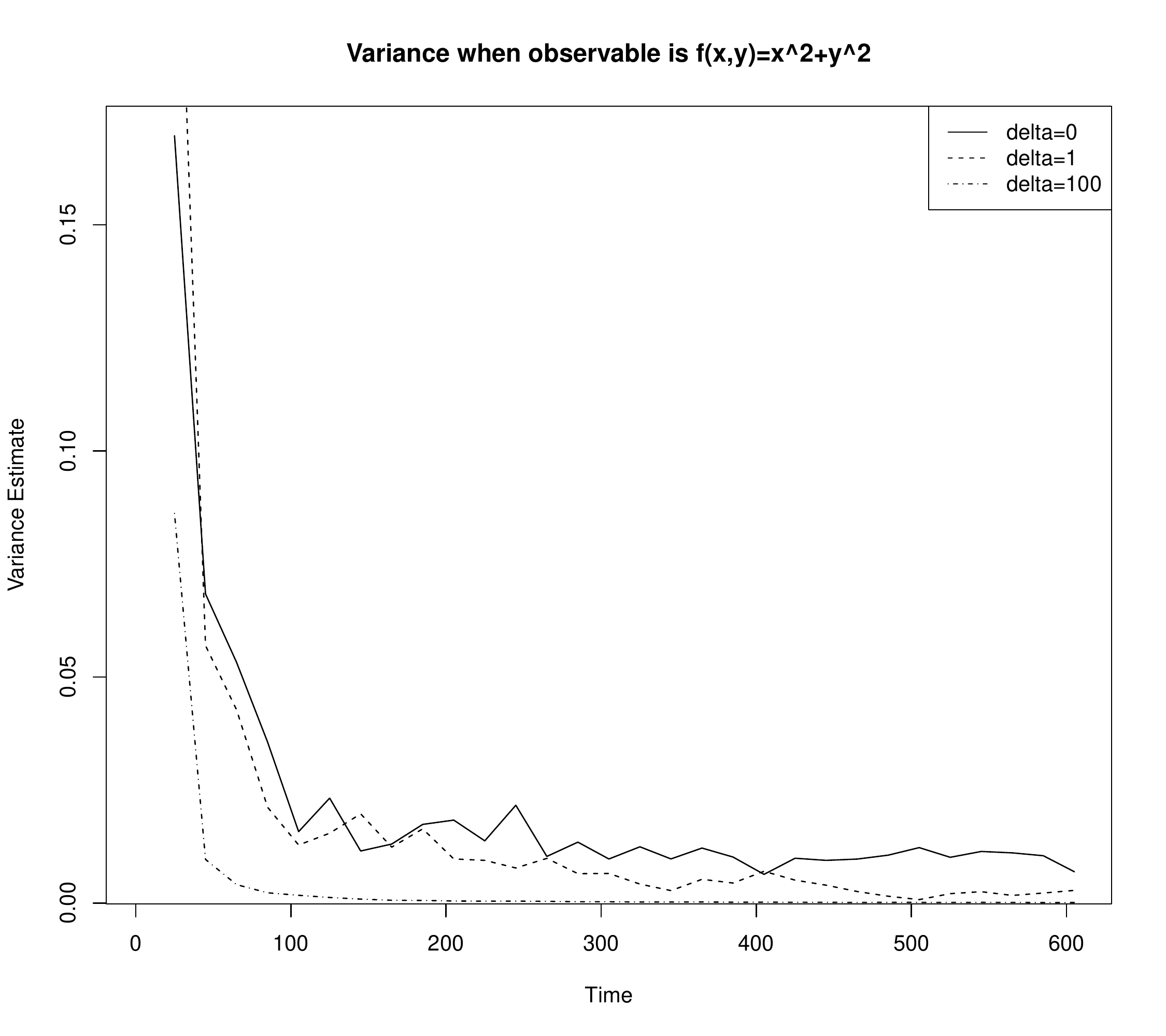}
\caption{Variance estimator for $f(x,y)=x^2 + y^2$} %
\label{Fig1a}%
\end{figure}

Lastly, to get a sense of the magnitude of the variance reduction, we take the observable $f(x,y)=x^2 + y^2$ (still $U(x,y)=\frac{1}{4}(x^{2}-1)^{2}+\frac{1}{2}y^{2}$) and
we use the standard batch means method, see \cite{AsmussenGlynn2007}, to estimate the asymptotic variance of $t^{-1} \int_0^t f(X_s) \, ds$, see Figure \ref{Fig1a}.
We present in Table \ref{Table1a1},
variance estimates for different values of $\delta=1/\epsilon $ and time horizon $t$. These numerical values are part of the values that were used to plot Figure \ref{Fig1a}.
It is clear that the variance reduction for this particular example is at the order of at least 2 magnitudes. Moreover, by Theorem \ref{T:LimitingVariance}, it is clear that for
small $\epsilon>0$, the variance estimate can be also considered as an approximation to the asymptotic variance of the corresponding estimation problem on the graph, see (\ref{Eq:LimitingVarianceOnTheGraph}).
\begin{table}[th]
\begin{center}
{\small
\begin{tabular}
[c]{|c|c|c|c|c|c|c|c|}\hline
$\delta\hspace{0.1cm} | \hspace{0.1cm} t$ & $25$ & $100$ & $200$ & $300$ & $400$ &
$500$ & $600$ \\\hline
$0$ &  $0.16$ & $0.015$ & $0.018$ & $0.010$ & $0.006$ & $0.012 $ & $0.007 $\\\hline
$1$ &  $0.25$ & $0.012$ & $0.009$ & $0.006$ & $0.005$ & $0.001 $ & $0.002 $\\\hline
$100$ &$0.09$ & $0.002$ & $4e-04$ & $2e-04$ & $2e-04$ & $1e-04 $& $1e-04$\\\hline
\end{tabular}
}
\end{center}
\caption{Estimated variance values for different pairs $(\delta,t)$.}%
\label{Table1a1}%
\end{table}



\begin{thebibliography}{0}

\bibitem{AsmussenGlynn2007} S. Asmussen and P.W. Glynn, {\em Stochastic Simulation}, Springer, 2007.


\bibitem{BrinFreidlin2000}
M. Brin and M. Freidlin. On stochastic behavior of perturbed Hamiltonian
systems, {\em Ergodic Theory and Dynamical Systems}, Vol.  20, (2000), pp. 55–76.

\bibitem{ConstantinKiselevRyshikZlatos2008}P. Constantin, A. Kiselev, L. Ryshik and A. Zlatos,
Diffusion and mixing in fluid flow
{\em Annals of Mathematics}, Vo.  168 (2008),  pp. 643-674


\bibitem{DonskerVaradhan1975}M.D. Donsker and S.R.S. Varadhan,  Asymptotic evaluation of certain Markov
process expectations for large times, I, {\em Communications Pure in Applied Mathematics}, Vol. 28,  (1975), pp. 1-47,
II, {\em Communications on Pure in Applied Mathematics}, Vol. 28,  (1975), pp. 279--301, and III, {\em Communications on Pure in Applied
Mathematics}, Vol. 29,  (1976), pp. 389-461.

\bibitem{DupuisDoll1}
P. Dupuis, Y. Liu, N. Plattner, and J. D. Doll, On the Infinite Swapping Limit for Parallel Tempering.
{\em SIAM Multiscale Modeling and Simulation}, Vol. 10, Issue 3, (2012), pp. 986-1022.

\bibitem{FrankeHwangPaiSheu2010}B. Franke, C.-R. Hwang, H.-M. Pai, and S.-J. Sheu, The behavior of the spectral gap under growing drift, {\em Transactions of the American Mathematical Society}, Vol 362, No. 3 (2010), pp. 1325-1350. 

\bibitem{FreidlinWeber2004}
M.I. Freidlin and  M. Weber, Random perturbations of dynamical systems and diffusion processes
with conservation laws, {\em Probability Theory and Related Fields}, Vol. 128, (2004), pp. 441-466.

\bibitem{FreidlinWentzell1993}
M. Freidlin and A.D. Wentzell. Diffusion Processes on Graphs and the
Averaging Principle, {\em Annals of Probability}, Vol. 21, No. 4, (1993), pp. 2215-2245.


\bibitem{FW2}
M.I. Freidlin and A.D. Wentzell, Random Perturbations of Hamiltonian
Systems, {\em Memoirs of the American Mathematical Society}, Vol. 109, No. 523, 1994.

\bibitem{Gartner1977}J. G\"{a}rtner, On large deviations from the invariant measure, {\em Theory of probability and its applications}, Vol. XXII, No. 1, (1977), pp. 24-39.

\bibitem{HwangMaSheu2005} C.-R. Hwang, S.Y. Hwang-Ma and S.-J. Sheu, Accelerating diffusions, {\em The Annals of Applied Probability}, Vol 15, No. 2, (2005), pp. 1433-1444.


\bibitem{HwangNormandWu2014} C.-R. Hwang, R. Normand and S.-J. Wu, Variance reduction for diffusions, {\em arXiv: 1406.4657}, 2014.


\bibitem{ReyBelletSpiliopoulos2014} L. Rey-Bellet and K. Spiliopoulos, Irreversible Langevin samplers and variance reduction: a large deviations approach,
{\em arXiv: 1405.0105}, 2014, submitted.
\end{thebibliography}
\end{document}